
\input amstex
\loadbold


\define\ap{Ann.\ Probab.}
\define\mprf{Markov Process.\ Related Fields}
\define\ptrf{Probab.\ Theory Related Fields}

\define\superzd{^{\raise1pt\hbox{$\scriptstyle {\Bbb Z}^d$}}}
\define\xsubvn{x_{\lower1pt\hbox{$\scriptstyle V(n)$}}}
\define\xssubvn{x_{\lower1pt\hbox{$\scriptscriptstyle V(n)$}}}

\define\ssub#1{_{\lower1pt\hbox{$\scriptstyle #1$}}} 
\define\ssubb#1{_{\lower2pt\hbox{$\scriptstyle #1$}}} 
\define\Ssub#1{_{\lower1pt\hbox{$ #1$}}} 
\define\Ssu#1{_{\hbox{$ #1$}}} 
\define\Ssubb#1{_{\lower2pt\hbox{$ #1$}}} 
\define\sssub#1{_{\lower1pt\hbox{$\scriptscriptstyle #1$}}} 
\define\sssubb#1{_{\lower2pt\hbox{$\scriptscriptstyle #1$}}} 
\define\sssu#1{_{\hbox{$\scriptscriptstyle #1$}}}                                              
                                     
\define\<{\langle}
\define\>{\rangle}

\define\({\left(}
\define\){\right)}
\define\[{\left[}                                   
\define\]{\right]}                                            
\define\lbrak{\bigl\{}
\define\rbrak{\bigr\}}
\define\lbrakk{\biggl\{}
\define\rbrakk{\biggr\}}


\define\osup#1{^{\raise1pt\hbox{$ #1 $}}}  
\define\osupm#1{^{\raise1pt\hbox{$\mskip 1mu #1 $}}}  
\define\osupp#1{^{\raise2pt\hbox{$ #1 $}}}  
\define\osuppp#1{^{\raise3pt\hbox{$ #1 $}}}  
\define\sosupp#1{^{\raise2pt\hbox{$\ssize #1 $}}}  
\define\sosup#1{^{\raise1pt\hbox{$\scriptstyle \mskip 1mu #1 $}}} 
\define\ssosup#1{^{\raise1pt\hbox{$\scriptscriptstyle #1 $}}} 
\define\ssosupp#1{^{\raise2pt\hbox{$\scriptscriptstyle #1 $}}} 
\define\ssosuppp#1{^{\raise3pt\hbox{$\scriptscriptstyle #1 $}}}


\define\Ga{\varGamma}


\define\Pitil{{\widetilde \Pi}}

\define\ztil{{\tilde z}}



\define\e{\varepsilon}

\define\epi#1{\text{epi}\,#1}

\define\pikkuhyppy{\vskip .1in}

\define\vinf{\|v_0\|_{\infty}}


\define\mmL{{\bold L}}

\define\mmN{{\bold N}}

\define\mmR{{\bold R}}

\define\mmZ{{\bold Z}}

\define\mmGa{{\boldsymbol \Gamma}}



\define\FF{{\Cal F}}



\documentstyle{amsppt}
\magnification=\magstep1
\document     
\baselineskip=12pt
\pageheight{43pc} 

\NoBlackBoxes

\centerline{\bf Perturbation of the Equilibrium for a Totally Asymmetric }

\pikkuhyppy

\centerline{\bf Stick Process in One Dimension }

\pikkuhyppy

$$\text{1999}$$

\hbox{}

\centerline{  Timo  Sepp\"al\"ainen
\footnote""{ Research partially supported by NSF grant DMS-9801085. }}
\hbox{}
\centerline{Department of Mathematics}
\centerline{Iowa State University}
\centerline{Ames, Iowa 50011, USA}
\centerline{seppalai\@iastate.edu}

\hbox{}

\hbox{}

\flushpar
{\it Summary.} We study the evolution of a 
small perturbation of the equilibrium of a
totally asymmetric one-dimensional interacting system. 
The model we take as example is Hammersley's process
as seen from a tagged particle, which can be 
viewed as a process of interacting positive-valued
stick heights on the sites of $\mmZ$. 
It is known that under Euler scaling (space
and time scale $n$) the empirical
stick profile obeys the Burgers equation.
 We refine this result in two ways: If the
process starts close enough to equilibrium, then
 over times $n^\nu$ for $1\le\nu<3$, and 
  up to errors that vanish in  hydrodynamic
scale, the dynamics merely translates
the initial stick configuration. 
A time evolution for the perturbation is visible  
under a particular family of scalings: 
 over times $n^\nu$, $1<\nu<3/2$,
 a perturbation of order $n^{1-\nu}$ 
from equilibrium follows the inviscid Burgers equation.
 The results for the stick model
are derived from asymptotic results for tagged particles
in Hammersley's process. 


\hfil

\hfil


\hfil

\flushpar

\hfil

\flushpar
Short Title: Perturbation of Equilibrium 

\break

\subhead 1.  Introduction \endsubhead
A number of recent papers have sought various refinements
to the basic hydrodynamic limits of 
interacting particle systems. One 
type of refinement is to add a small perturbation
to the equilibrium of the process, and study the 
time evolution of the perturbation. For asymmetric
exclusion in dimensions 3 and higher, Esposito, Marra and Yau (1994)
proved that under diffusive scaling 
(time scale $n^2$ and space scale $n$)
 a perturbation of order $n^{-1}$ follows a
conservation law with a diffusion term. The backdrop 
of this result is the standard hydrodynamic limit
of asymmetric processes, which leads to a 
conservation law {\it without} diffusion term under {\it Euler
scaling} (time scale and space scale both $n$). 
A context for the Esposito et al.\ result is  
the search for microscopic interpretations of 
  the Navier-Stokes 
equations. We refer the 
reader to p.\ 185--188 in Kipnis and Landim's (1999)
monograph for a description of this program and further references.

Our paper looks at the question of 
Esposito  et al.\ in one dimension. We
add a perturbation of order $n^{-\beta}$ to the 
equilibrium, $\beta>0$.
 The perturbation 
 vanishes in the hydrodynamic limit $n\to\infty$,
and we study the effect of this perturbation 
under various  time scales $n^\nu t$, $\nu\ge 1$. We have two types of 
results: (1) For $\beta\in(0,1/2)$,
a hydrodynamic limit in the 
time scale $n^{1+\beta}t$ shows that the 
 perturbation obeys macroscopically the Burgers  
equation {\it without} diffusion term. (2) For 
$\nu\in[1,3)$ and $\beta$ close enough to 1, 
we show that the dynamics is simply a translation of the 
initial configuration, up to $o(n)$ error terms.

The most popular models  for studies of the hydrodynamics 
of asymmetric stochastic 
dynamics are the exclusion and zero-range
processes. Instead of these processes, we prove our results
for the so-called {\it stick process}, which can also be
regarded  as Hammersley's process as seen from a 
tagged particle. This process 
has nonnegative variables 
$\eta(i)$ [stick heights] on the sites of $\mmZ$ that 
exchange pieces between each other. The stick process lacks 
some of the good properties of the 
exclusion or zero-range process: the state space
is not compact, the rates are unbounded, 
and the amount of material that jumps
is also unbounded. 

To make up for these 
complications the totally asymmetric one-dimensional 
process has a beautiful combinatorial structure
uncovered by Aldous and Diaconis (1995). 
This structure connects Hammersley's process
and the stick model to the increasing sequences 
problem on planar Poisson points. A key
ingredient of our proofs are sharp deviation 
estimates for the inceasing sequences problem 
 from Kim (1996), 
Sepp\"al\"ainen (1998b), and Baik, Deift and 
Johansson (1999). 

We believe that the results
of our paper hold also for totally
asymmetric exclusion and zero-range processes. 
The basis for this conjecture is that these
processes possess particle-level variational formulations
that involve planar growth models, 
analogous to the increasing sequences
 connection of Hammersley's 
process [Sepp\"al\"ainen (1998a,c)]. 
Johansson (1999) has shown that the limiting 
fluctuations for this growth model are the same 
as for the increasing sequences model.

{\it Organization of the paper. } In section 2 we describe 
the stick model and state the results mentioned 
above: Theorem 1 gives
the translation, and Theorem 2 the hydrodynamic limit
of the perturbation. Theorems 1 and 2 are corollaries
of corresponding Theorems 3 and 4 for tagged particles in 
Hammersley's process. These are stated in Section 3. 
The translation Theorem 3 for Hammersley's process
is compared to a similar result of Ferrari and Fontes
(1994) for asymmetric exclusion. 
Section 4 addresses briefly the rigorous construction 
of Hammersley's process and the stick process, and 
the connection with increasing sequences. Sections 5--7
contain the proofs: Section 5 contains lemmas, 
Section 6 the proof of Theorem 3, and 
Section 7 the proof of Theorem 4.

\hbox{}

\subhead 2. The stick model and the results \endsubhead
Here is an informal description of the model. 
A rigorous construction will follow in Section 4. 
The state of the process
 is a configuration $\eta=(\eta(i):i\in\mmZ)$
where each $\eta(i)$ is a nonnegative real number.
Think of $\eta(i)$ as the height of a vertical stick
attached to site $i\in\mmZ$. At exponential rate 
equal to $\eta(i)$, the following event takes place:
Pick a random quantity $u$ uniformly distributed 
on $[0,\eta(i)]$. Break off a piece of length $u$ 
from the stick at $i$, and attach this piece to
the stick at site $i+1$. Thus if the neighboring stick
lengths before the event were $(\eta(i), \eta(i+1))$,
then after the event they are  $(\eta(i)-u,\eta(i+1)+u)$. 
These events happen at all sites $i$ independently of
each other.  In the language of generators, this
dynamics is  expressed as
$$Lf(\eta)=\sum_{i\in\mmZ}\int_0^{\eta(i)}
[f(\eta^{u,i,i+1})-f(\eta)]du\,,
\tag 2.1
$$
where $\eta^{u,i,i+1}=\eta-u\delta_i+u\delta_{i+1}$. 

In Sepp\"al\"ainen (1996) a Markov
 process
$\eta(t)=(\eta(i,t):i\in\mmZ)$, $t\ge 0$, is
constructed that operates according to this description. 
The state space of the process is 
$$Y=\lbrakk \eta\in[0,\infty)^\mmZ:
\lim_{N\to-\infty}N^{-2}\sum_{i=N}^{-1}\eta(i)=0\,\rbrakk
\tag 2.2
$$
and the paths of the process are in the Skorohod space
$D([0,\infty),Y)$.
$Y$ is not closed in the product topology, but 
is given a stronger topology with a complete,
separable metric. $L$ of (2.1) is the generator of 
the process, in the sense that
$$E^\eta\bigl[f(\eta(t))\bigr]-f(\eta)=\int_0^t 
E^\eta\bigl[Lf(\eta(s))\bigr]ds
\tag 2.3
$$
 for all bounded continuous cylinder 
functions $f$ on $Y$ and all initial states $\eta\in Y$. 
 $E^\eta$ stands for the expectation under the 
path measure of the process started at state $\eta$. 
Furthermore, the process has a 1-parameter family of
invariant distributions, namely the i.i.d.\ exponential
distributions on the variables $(\eta(i):i\in\mmZ)$. 

We focus now on the hydrodynamic behavior of this process.
The basic result  [Sepp\"al\"ainen (1996)] is that 
under Euler scaling the empirical stick profile 
obeys the Burgers equation. Suppose $u(x,t)$, 
$(x,t)\in\mmR\times[0,\infty)$, is an entropy solution
of the Burgers equation
$$u_t+ (u^2)_x=0\,,
\quad u(x,0)=u_0(x)
\tag 2.4
$$
with nonnegative initial data $u_0\in L^\infty(\mmR)$. 
Consider a sequence $\eta^n$, $n=1,2,3,\ldots\,$, of
stick processes, and assume  that a law of large numbers 
is satisfied at time $t=0$:
for all $a<b$ in $\mmR$, 
$$\lim_{n\to\infty}  n^{-1}\sum_{i=[na]+1}^{[nb]}\eta^n(i,0)
=\int_a^b u_0(x)dx\qquad\text{in probability.}
\tag 2.5
$$
The theorem is that the law of large numbers 
continues to hold at  all later times $t>0$: 
$$\lim_{n\to\infty}  n^{-1}\sum_{i=[na]+1}^{[nb]}\eta^n(i,nt)
=\int_a^b u(x,t)dx\qquad\text{in probability.}
\tag 2.6
$$

Euler scaling refers to the scaling in the above limit, 
where the ratio of macrosopic and microscopic units is the
same $n$ for both space and time.
A macroscopic space interval $(a,b]$ corresponds
to approximately $n(b-a)$ microscopic lattice sites, and macrosopic
time $t$ corresponds to microscopic time $nt$. 
The derivation of the hydrodynamic limit (2.6) from 
the hypothesis (2.5) requires some technical assumptions,
and the details can be found in Sepp\"al\"ainen (1996).

A trivial special case of the hydrodynamic limit 
is of course the case of a process in equilibrium:
If the sticks are initially i.i.d.\ exponentially
distributed with common expectation $E[\eta(i,0)]=q$,
then this situation persists, and the 
macroscopic profile is the constant $u(x,t)\equiv q$. 

In the present paper we study the evolution
of a small perturbation 
of the equilibrium. 
The initial macroscopic profile is 
$$u_0(x)=q+n^{-\beta}v_0(x)
\tag 2.7
$$
where $q>0$ is the fixed equilibrium density, 
$v_0$ is a bounded measurable function on $\mmR$,
and $\beta\in(0,1]$ is a parameter that we adjust
to investigate different scalings. The function 
$v_0$ is not assumed to take any particular sign, 
so to have nonnegative profiles
 we consider only
$n$ large enough to have $q>n^{-\beta}\vinf$. 
For each $n$, the 
initial stick configuration $(\eta^n(i,0):i\in\mmZ)$ is assumed to be in 
{\it local equilibrium} with 
macroscopic profile $u_0$. Precisely speaking, our
assumption is this: 
$$
\aligned
&\text{The 
variables $(\eta^n(i,0):i\in\mmZ)$ are mutually independent,}\\
&\text{exponentially distributed, and have expectations}\\
&E\bigl[ \eta^n(i,0)\bigr]=q+n^{1-\beta}\int_{(i-1)/n}^{i/n}v_0(x)dx\,.
\endaligned
\tag 2.8
$$
The perturbation of the expected density is taken to be
$$
\text{$n^{-\beta}\cdot\{$the average of $v_0$ over the interval
$((i-1)/n,i/n]\,\}$,}
$$
 instead of $n^{-\beta}\cdot\{$the point value
$v_0(i/n)\}$ because we are not making any regularity  assumptions
on $v_0$. A standing assumption is also that $q>0$. 
In Section 4.1 we explain how the 
case $q=0$ is reduced to the basic hydrodynamic limit
(2.5)--(2.6). 

Since $\|u_0-q\|_\infty\longrightarrow 0$ as $n\to\infty$, 
 the limit (2.6) is valid again with constant
profile $u(x,t)\equiv q$.  To
escape the regime of (2.6) we subtract the equilibrium 
density $q$ and 
 speed up time more, beyond the hydrodynamic scale 
$nt$. We introduce a second parameter $\nu\in[1,\infty)$, 
and look at the evolution of the stick profile over times 
of order $n^\nu$. The space scaling will be the same 
as in (2.5)--(2.6), so the lattice of sites 
scales as $n^{-1}\mmZ$. 
Our object of study is the empirical 
 stick profile 
$$\sum_{i\in\mmZ}\lbrak\eta^n(i,n^\nu t)-q\rbrak \delta_{i/n}\,.$$
In other words, we follow either integrals 
$\sum_{i\in\mmZ}\lbrak\eta^n(i,n^\nu t)-q\rbrak\phi(i/n)$
of  compactly supported, continuous test functions $\phi$,
or equivalently the total stick mass in  macroscopic intervals
$(x,y]$, $\sum_{i=[nx]+1}^{[ny]}\eta^n(i,n^\nu t)-nq(y-x)$. 

Let us derive an easy ``benchmark'' result against
which we can compare later results. 
It is proved in Sepp\"al\"ainen (1996) that the 
stick process is {\it attractive}. This means that,
if $\eta$ and $\zeta$ are two initial states that 
satisfy $\eta\ge\zeta$ [inequalities are interpreted
coordinatewise,  
$\eta(i)\ge \zeta(i)$ for all $i\in\mmZ$], then it is possible
to construct the processes $\eta(t)$ and $\zeta(t)$ 
on a common probability space so that the inequality
$\eta(t)\ge \zeta(t)$ holds at all times $t\ge 0$, almost
surely. 

Fix $n$ for the moment. 
Let $\zeta^1$ and $\zeta^2$ be stick processes in
equilibrium, with expectations 
$$
\text{$E[\zeta^1(i,t)]=q-n^{-\beta}\vinf\ $
and $\ E[\zeta^2(i,t)]=q+n^{-\beta}\vinf$.}
$$
In other words, for each fixed time $t$
and for $r\in\{1,2\}$, the stick heights
$(\zeta^r(i,t):i\in\mmZ)$ are  exponentially
distributed i.i.d.\ random variables with expectations 
as above. At time $t=0$, 
 we can construct the initial configurations
of all three processes  $\zeta^1$, $\eta^n$, and $\zeta^2$
on a single probability
space so that
$$\zeta^1(i,0)\le\eta^n(i,0)\le\zeta^2(i,0)\qquad\text{for all $i$,
a.s.}$$ 
To do this take an i.i.d.\ sequence of Exp(1) variables $\{X_i\}$ and 
set 
$$\text{
$\zeta^r(i,0)=E[\zeta^r(i,0)] X_i$ and 
$\eta^n(i,0)=E[\eta^n(i,t)] X_i$ for $r=1,2$ and all $i$.}
$$
We then construct all three processes 
on a common probability space so that
$\zeta^1(i,t)\le\eta^n(i,t)\le\zeta^2(i,t)$ for all $t$ and $i$,
 a.s. 

Let $\{K_n\}$ be an
arbitrary sequence of integers, to be used as 
translations on the lattice. 
The construction gives these inequalities:
$$\aligned
&\sum_{i=[nx]+1}^{[ny]} \bigl\{ \zeta^1(K_n+i,n^\nu t)
-E\bigl[\zeta^1(i,n^\nu t)\bigr]\bigr\} 
-([ny]-[nx])n^{-\beta}\vinf\\
&\le \sum_{i=[nx]+1}^{[ny]} \bigl\{ \eta^n(K_n+i,n^\nu t)
-q\bigr\} \\
&\le \sum_{i=[nx]+1}^{[ny]} \bigl\{ \zeta^2(K_n+i,n^\nu t)
-E\bigl[\zeta^2(i,n^\nu t)\bigr]\bigr\} 
+([ny]-[nx])n^{-\beta}\vinf\,.
\endaligned
\tag 2.9
$$
A terminological convention:
Throughout, we shall say 
$$X=Y+o(n^\alpha)\qquad\text{a.s.}
\tag 2.10
$$
as a shorthand for 
$$\lim_{n\to\infty} n^{-\alpha}\bigl|X-Y\bigr|=0\qquad\text{a.s.}
\tag 2.11
$$
The sums of the  $\zeta$-terms in (2.9)
are almost surely $o(n^{1/2+\delta})$ for any $\delta>0$
because the variables  are i.i.d. and have sufficient moments. 

Thus we get the result
$$ \sum_{i=[nx]+1}^{[ny]} \eta^n(K_n+i,n^\nu t)
= nq(y-x)+o\bigl(n^{[1/2]\vee[1-\beta]+\delta}\bigr)
\qquad\text{a.s.}
\tag 2.12
$$
for any $\nu$, $\beta>0$, and
arbitrarily small $\delta>0$. The 
translation $K_n$ was included in anticipation
of later results.  Because
we are speeding up time beyond the hydrodynamic scale,
a certain translation will appear naturally.
 The goal of the remainder of the paper is to
improve on (2.12), by obtaining results that reveal how the
perturbation evolves in time, or have a smaller error term.

As the last preparatory step,
we  construct the solution of the 
Burgers equation (2.4) by the Hopf-Lax formula. 
The perturbation  $v_0(x)$ is now the
initial data. Define $V_0(x)$ by 
$$
V_0(0)=0\quad\text{ and }\quad V_0(y)-V_0(x)=\int_x^y v_0(z)dz
\quad\text{for all $x<y$.}
$$
$V_0$ is a Lipschitz function with a bounded derivative a.e.
For $(x,t)\in\mmR\times[0,\infty)$, define 
$V(x,0)=V_0(x)$ and for $t>0$ 
$$V(x,t)=\inf_{y\in\mmR}\lbrakk V_0(y)
+\frac1{4t}{(x-y)^2}\rbrakk\,.
\tag 2.13
$$
Then $V$ is the unique  {\it viscosity solution} of the 
Hamilton-Jacobi equation 
$$V_t+(V_x)^2=0\,,\quad V(x,0)=V_0(x)\,.
\tag 2.14
$$
For each fixed $t$, $V(\cdot\,,t)$ is again a Lipschitz
function, so it has a.e.\ an $x$-derivative $v=V_x$. 
This function $v(x,t)$ is the unique {\it entropy solution}
of (2.4) with initial data $v_0$. The reader can find
a development of these p.d.e.\ results
 in Evans (1998). 

Now the results for the stick process. 
 The most general 
result, valid for all scalings, 
does not  identify any 
time evolution, only a translation of the initial 
sticks. 

\proclaim{Theorem 1} Assume that  $\beta>0$ and $\nu\ge 1$. 
Let $\eta^n(t)$ denote
the stick process started from the
initial configuration {\rm (2.8)}. Fix
$x<y$ in $\mmR$ and $t>0$. Then, for any $\delta>0$, 
 the following asymptotic
equality is valid almost surely as $n\to\infty$:  
$$ \sum_{i=[nx]+1}^{[ny]}  \eta^n([2n^\nu qt]+i,n^\nu t)
 = \sum_{i=[nx]+1}^{[ny]}\eta^n(i,0) 
+o(n^{[\nu-2\beta]\vee[\nu/3]+\delta})\,.
\tag 2.15
$$
\endproclaim

Why the translation $[2n^\nu qt]$ appears 
naturally is explained in Section 4.2. 
The error exponent in the statement (2.15) satisfies
$$[\nu-2\beta]\vee[\nu/3]=\cases
\nu-2\beta, &\nu>3\beta\\
\nu/3, &\nu\le 3\beta\,.
\endcases
\tag 2.16
$$
If $\nu<3$ and $\beta>(\nu-1)/2$ the error in (2.15) is
$o(n)$ and so vanishes in the standard hydrodynamic scaling of (2.6). 

We find one family of scalings where the perturbation 
evolves according to the Burgers equation. For this to
happen the perturbation has to be larger than $n^{-1/2}$. 

\proclaim{Theorem 2} Suppose $\beta\in(0,1/2)$
and set $\nu=1+\beta$. Let $\eta^n(t)$ denote
the stick process started from the
initial configuration {\rm (2.8)}.
Let $\phi$ be a compactly
supported, continuous test function on $\mmR$. 
Then almost surely
$$\lim_{n\to\infty}\frac1{n^{1-\beta}}\sum_{i\in\mmZ}
\bigl\{ \eta^n([2n^{1+\beta} qt]+i,n^{1+\beta} t)
-q\bigr\} \phi(i/n)=\int_{\mmR}\phi(x)v(x,t)dx\,. 
\tag 2.17
$$
\endproclaim

This result compares directly with Corollary 2.3 in
Esposito et al.\ (1994),
where the corresponding result is proved 
in dimensions $d\ge 3$ for an exclusion process. 
The deterministic limit
(2.17) cannot be valid for $\beta=1/2$ 
because in equilibrium this would be the 
 central limit theorem scaling. 

\demo{ Remark about construction}
The almost sure results of Theorems 1 and 2  are proved 
for a special construction explained in Section 4. 
In this construction the processes $\eta^n$ 
are defined on one common probability space, and 
 the variables $\eta^n(i,t)$ 
are realized as  interparticle distances
of Hammersley's process. 
This is {\it not} the construction used in (2.9)
that makes $\eta$ 
attractive. Both theorems are  proved by  Borel-Cantelli
arguments, and the probability estimates 
 for the arguments are derived  with the help
of the special
construction. But once derived, 
the estimates are valid in all
constructions because they are statements about
the  distributions of the processes.
 Hence Theorems 1 and 2 are valid for any construction
of the stick process. 
\enddemo

To compare (2.17) directly with (2.12), we can write it
in the form 
$$ \aligned
&\sum_{i=[nx]+1}^{[ny]} \eta^n([2n^{1+\beta} qt]+i,n^{1+\beta} t)\\
&\qquad\qquad= nq(y-x)+n^{1-\beta} \int_x^y v(z,t)dz
  +o\bigl(n^{1-\beta}\bigr)\,.
\endaligned
\tag 2.18
$$
A comparison of the errors in (2.12), (2.15), and (2.18) reveals
that for $\nu>1+\beta$ the easy result (2.12)  in fact has
 the 
smallest error. For $\nu=1+\beta$ (2.18) has the smallest error, 
while for $\nu<1+\beta$ it can be either one of (2.12) and
 (2.15), depending
on the exact relation between $\beta$ and $\nu$. 

\hbox{}

\subhead 3.\ Asymptotics for a tagged particle in Hammersley's process
\endsubhead
In Hammersley's process a countable collection of point
particles evolves on $\mmR$ according to the following
 rule: if $x<y$ are two locations of neighboring
particles, then with rate equal to the distance $y-x$
the particle at $y$ jumps to a randomly (uniformly) 
chosen location in the interval $(x,y)$. All particles
execute jumps independently of each other. 

This 
evolution can be graphically constructed 
with a rate one homogeneous Poisson point process
on the space-time  plane $\mmR\times(0,\infty)$:
Suppose $(x,t)$ is a point of the Poisson process. 
Then  at time $t$
the leftmost particle in  $[x,\infty)$
jumps to $x$. If the leftmost particle already
happened to be at $x$, or if there is no leftmost particle 
in $[x,\infty)$, no jump takes place. This
latter case can happen
if there are infinitely many particles in some bounded
interval. 

There is an obvious connection between Hammersley's process
and our stick process.  We  assume that we
can label the particles by integers
in an order-preserving way. Let $z(i,t)$ denote the 
position of particle $i$ at time $t$. The assumption is 
$$z(i-1,t)\le z(i,t)\qquad\text{for all $i$ and $t$.}
\tag 3.1
$$ 
Suppose we have constructed 
 the process $z(t)=(z(i,t):i\in\mmZ)$
that operates  according to the
description above. 
Define 
$$\eta(i,t)=z(i,t)-z(i-1,t)\qquad\text{for $i\in\mmZ$}\,. 
\tag 3.2
$$
Then it is clear that $\eta(t)$ evolves as our stick 
process. When particle $z(i)$ jumps to the left,
stick $\eta(i)$ donates a piece to stick $\eta(i+1)$. 
In particle system jargon, the stick process is 
Hammersley's process as seen from a tagged particle. 
What this means is that knowing $\eta(t)$
and the evolution of one particle $z(j,t)$ is equivalent
to knowing the process $z(t)$. 
 The simultaneous construction of Hammersley's process
and the stick process is discussed
in Section 4. 

Assume that
 the initial sticks $\eta^n(0)=(\eta^n(i,0):i\in\mmZ)$
that satisfy (2.8) have been defined on some probability space.
Initial particle configurations 
$z^n(0)=(z^n(i,0):i\in\mmZ)$ are defined
on this same probability space  by 
$$\aligned
&z^n(0,0)=0\,,\,\,z^n(i,0)=\sum_{j=1}^i\eta^n(j,0)\quad\text{ for $i>0$},\\
 &\text{and }\quad z^n(i,0)=\sum_{j=i+1}^0\eta^n(j,0)\quad\text{ for $i<0$.}
\endaligned 
\tag 3.3
$$ 
Thus $z^n(i,0)$ is a sum of independent exponential 
random variables with uniformly bounded expectations, 
and 
$$E\bigl[z^n(i,0)\bigr]=qi+n^{1-\beta}V_0(i/n)\,.
\tag 3.4
$$

The processes $\{z^n(t)\}$ are then constructed together on one
probability space where the initial configurations
$\{z^n(0)\}$ and the space-time Poisson points are independent. 
All processes $z^n(t)$ use
 the same realization of the Poisson points to construct 
the dynamics. This is not really necessary because our
a.s.\ results come from Borel-Cantelli arguments. But in 
the proof it is convenient to work with a single Poisson
process and the family $\{z^n(0)\}$ of initial configurations, 
instead of giving each process $z^n(t)$ 
its own space-time Poisson process.

\proclaim{Theorem 3} Assume that  $\beta>0$ and $\nu\ge 1$. 
Let $z^n(i,t)$ denote
 Hammersley's process started from the
initial configuration described in {\rm (3.3)} and {\rm (2.8)}. Fix
$x\in\mmR$ and $t>0$. Then, for any $\delta>0$,  
 we have the following asymptotic
equality almost surely as $n\to\infty$: 
$$z^n([nx]+[2n^\nu qt],n^\nu t)=n^\nu tq^2 +z^n([nx],0) 
+o(n^{[\nu-2\beta]\vee[\nu/3]+\delta})\,.
\tag 3.5
$$
\endproclaim

Theorem 1 is an immediate consequence of (3.2) and  Theorem 3. 

Ferrari and Fontes (1994) proved a translation
result of this type for the exclusion process. Suppose
for the moment that  the 
$\eta(i)$'s are occupation variables of
totally asymmetric 1-dimensional
simple exclusion in equilibrium at density $\rho$. 
Then the jumps of the  $z$-variables correspond to the current 
of particles. 
Statement (1.5) in Theorem 1 of Ferrari  and Fontes (1994)
implies that, in the  $L^2$ sense as $n\to\infty$, 
$$z(0,nt)=nt\rho^2 +z([nth(\rho)],0)+o(n^{1/2})\,,
\tag 3.6
$$
where $h(\rho)=2\rho-1$. 
This can be compared with our result 
for Hammersley's process: 
with $\nu=1$ and 
 $\beta\ge 1/3$, (3.5)
 implies  that 
$$z^n(0,n t)=n tq^2 +z^n(-[2tqn],0) 
+o(n^{1/3+\delta})\,.
\tag 3.7
$$
The error is smaller in (3.7)  than in (3.6), but the
Ferrari-Fontes result is valid for more general
asymmetric exclusions,  not only for totally
asymmetric.  

Next a result with explicit time evolution. 
In one of the cases treated by the next theorem, we 
will assume that $V_0(x)$ 
has asymptotic slopes in the sense that these limits exist:
$$v_0(-\infty)=\lim_{x\to-\infty}\frac{V_0(x)}{x}
\quad\text{and}\quad
v_0(+\infty)=\lim_{x\to+\infty}\frac{V_0(x)}{x}\,.
\tag 3.8
$$
When this is the case, we define the piecewise linear
``asymptotic profile''
$$V_\infty(x,0)=\cases 
v_0(-\infty)x\,, &x<0\\
0\,,&x=0\\
v_0(+\infty)x\,, &x>0\\
\endcases
\tag 3.9
$$
and its evolution for $t>0$ by
$$V_\infty(x,t)=\inf_{y\in\mmR}\lbrakk V_\infty(y,0)
+\frac1{4t}{(x-y)^2}\rbrakk\,.
\tag 3.10
$$

\proclaim{Theorem 4} Assume $\nu>3\beta$, in addition to the basic
assumption $\beta>0$ and $\nu\ge 1$. Let $z^n(i,t)$ denote
 Hammersley's process started from the
initial configuration $z^n(i,0)$ described in {\rm (3.3)} and 
{\rm (2.8)}. Fix
$x\in\mmR$ and $t>0$. Then we have the following asymptotic
equalities, each statement valid almost surely as $n\to\infty$. 

Case 1: $\nu>1+\beta$. Assume that the limits in {\rm (3.8)} exist.  Then 
$$z^n([nx]+[2n^\nu qt],n^\nu t)=n^\nu tq^2 +nxq+n^{\nu-2\beta}V_\infty(0,t)
+o(n^{\nu-2\beta})\,.
\tag 3.11
$$

Case 2: $\nu=1+\beta$. Then 
$$z^n([nx]+[2n^\nu qt],n^\nu t)=n^\nu tq^2 +nxq+n^{1-\beta}V(x,t)
+o(n^{1-\beta})\,.
\tag 3.12
$$

Case 3: $1\le \nu<1+\beta$. Then for any $\delta>0$, 
$$z^n([nx]+[2n^\nu qt],n^\nu t)=n^\nu tq^2 +nxq+n^{1-\beta}V_0(x)
+o(n^{[1/2]\vee[\nu-2\beta]+\delta})\,.
\tag 3.13
$$
\endproclaim

\demo{Remarks} 
Recall again (2.10)--(2.11) for the precise meaning of the almost sure 
$o(n^\alpha)$ error terms. 
In Case 1, the term $nxq$ in (3.11) may or may not be included in the 
error $o(n^{\nu-2\beta})$, depending on whether $\nu>1+2\beta$ or not. 
The statement (3.11) for Case 1 does not improve (2.12) because the 
error $n^{\nu-2\beta}$ is strictly larger than $n^{1-\beta}$
in this case. Theorem 2 follows from Case 2 above.

The remark about construction at end of Section 2 applies 
here too. The proofs of Theorems 3 and  
4 are Borel-Cantelli arguments that depend on
estimates of the distributions of the processes,
and hence are valid in all constructions. 
\enddemo

The three cases reveal the effect of the time scale 
on the evolution of the perturbation: For fast times 
$\nu>1+\beta$ we only see the asymptotic effect 
$V_\infty(0,t)$ which is independent of the reference 
point $x$. For slow times $\nu<1+\beta$ we only see the 
initial perturbation $V_0(x)$. And
 exactly at $\nu=1+\beta$, we see the perturbation 
evolve according to the Burgers equation.

It remains to prove Theorems 3 and 4. This proof uses
on a special construction of Hammersley's process in terms
of increasing sequences of the space-time Poisson points. 

\hbox{}

\subhead 4.\ Graphical construction and increasing sequences\endsubhead
Consider a planar, rate one, homogeneous
 Poisson
point process. A sequence $(x_1,t_1)$, $(x_2,t_2)$, $\ldots$, 
$(x_m,t_m)$ of Poisson points is {\it increasing} if 
$$x_1<x_2<\cdots<x_m\qquad\text{and}\qquad t_1<t_2<\cdots<t_m\,.
$$
For arbitrary $(a,s)$, $(b,t)$ on the plane, define the 
random variable $\mmL\bigl((a,s), (b,t)\bigr)$ as the 
maximal number of  Poisson points on an increasing sequence
contained in
the rectangle $(a,b]\times(s,t]$. Abbreviate
$\mmL(b,t)=\mmL\bigl((0,0), (b,t)\bigr)$ for the case
where the lower left corner is the origin. 

 An inverse to $\mmL$ is defined by
$$\mmGa\bigl((a,s), m, \tau\bigr)=\inf\bigl\{ h>0: 
\mmL\bigl((a,s), (a+h,s+\tau)\bigr) \ge m\bigr\}\,.
\tag 4.1
$$
In words: $\mmGa\bigl((a,s), m, \tau\bigr)$ is the 
minimal horizontal distance $h$ for which
the rectangle $(a,a+h]\times(s,s+\tau]$ contains
 an increasing
sequence of $m$ points. 
Again abbreviate 
$\mmGa(m,\tau)=\mmGa\bigl((0,0), m,\tau\bigr)$. 

These random variables  satisfy laws of large numbers: 
$$\lim_{s\to\infty} 
\frac1s\mmL(sb,st)= 2\sqrt{bt\,}\qquad\text{ and }\qquad
\lim_{s\to\infty} \frac1s\mmGa([sa],st)= \frac{a^2}{4t}\qquad\text{a.s.}
\tag 4.2
$$ 
The 
existence of the limits 
follows from the subadditive ergodic theorem.
The exact values were 
first calculated by Vershik and Kerov (1977). 

In the previous section we suggested how to construct 
Hammersley's process with a rate one space-time Poisson
point process. The rule was that Poisson point $(x,t)$ 
pulls the leftmost particle in  $[x,\infty)$ to the 
location $x$ at time $t$. 
As usual in particle system contexts, 
constructing the process rigorously from this
description, on the infinite real line,  needs a proof. 

We can take an elegant way out with the help of
the increasing paths. Assume given an initial 
configuration $z(0)=(z(i,0):i\in\mmZ)$ that satisfies
the ordering convention (3.1) 
for $t=0$. Given a realization of the Poisson points, 
{\it define }
$$z(k,t)=\inf_{i:i\le k}\lbrakk z(i,0)+
\mmGa\bigl((z(i,0) ,0), k-i, t\bigr)\rbrakk
\tag 4.3
$$
for all $k\in\mmZ$ and $t>0$. In words: The potential 
locations of $z(k,t)$ are all points $x$ such that the
rectangle $(z(i,0),x]\times(0,t]$ contains an increasing
sequence of $k-i$ Poisson points. Of these potential 
locations $z(k,t)$ chooses the leftmost. 

If we permit $-\infty$ 
as a value for $z(k,t)$, (4.3) defines a process 
$z(t)=(z(k,t):k\in\mmZ)$ that satisfies (3.1). 
To rule out the possibility of jumping to $-\infty$ in finite 
time, define the state space 
$$Z=\lbrak z=(z(i))\in\mmR^{\mmZ}: \text{$z(i-1)\le z(i)$
for all $i$, and }
\lim_{i\to-\infty} i^{-2}z(i)=0\rbrak\,.
\tag 4.4
$$
One can check that if $z(0)\in Z$, then almost surely 
the infimum in (4.3) is always attained at some finite $i$
and $z(t)\in Z$ for all
$t$. 
Homogeneity of the space-time Poisson point process then 
implies that (4.3) defines a time-homogeneous Markov process $z(t)$
with state space $Z$. 

Definitions (2.2) and (4.4) show that 
 $Y$ maps injectively into $Z$ through equations (3.3), and 
$Z$ back onto $Y$  through equation (3.2).  
So given an initial stick configuration $(\eta(i,0):i\in\mmZ)$
in $Y$, we define an initial particle configuration 
$z(0)\in Z$ as in (3.3), then define the process $z(t)$ by (4.3), 
and finally  use (3.2) to define 
the stick process $\eta(t)$. 

This is convenient as a rigorous definition of the 
processes $\eta(t)$ and $z(t)$, but 
it is not so obvious that the  resulting
dynamics  follows our earlier descriptions. 
One can prove that when $\eta(t)$ is defined 
 this way, (2.3) is satisfied so the generator of 
$\eta(t)$ is $L$. All the facts mentioned here can be found
in Sections 3--5 in Sepp\"al\"ainen (1996).    

We can also argue from definition (4.3)  that if
$(x,t)$ is a space-time Poisson point, then at time
$t$ the leftmost particle in $[x,\infty)$
is at $x$, if such a particle exists.  Suppose
not, so that for some $k$, 
$z(k-1,t)<x<z(k,t)$. Pick $i\le k-1$ so that 
$$z(k-1,t)=z(i,0)+
\mmGa\bigl((z(i,0) ,0), k-1-i, t\bigr)\,.
$$
 Barring the null event that
space-time Poisson points can lie on the same 
horizontal line, there must be an increasing 
sequence of $k-1-i$ Poisson points from $(z(i,0) ,0)$ to 
a point $(y,s)$ such that $y=z(k-1,t)<x$ and $s<t$. 
[In the extreme case $i=k-1$, this sequence is empty,
 and $s=0$.] 
Consequently we can append the new point $(x,t)$ 
to this increasing sequence to produce a sequence
of $k-i$ points from $(z(i,0) ,0)$ to $(x,t)$.
Then 
$$z(k,t)\le z(i,0)+
\mmGa\bigl((z(i,0) ,0), k-i, t\bigr)\le x\,,
$$
contradicting $x<z(k,t)$.

The remainder of the paper proves Theorems 3 and 4 
through definition (4.3). The construction
of the family of processes $\{z^n(t)\}$ is the 
following. There is a single probability space $(\Omega,\FF,P)$
on which are defined the initial
locations $\{z^n(0)\}$ and, independently of them, 
 the space-time Poisson
point process.
On this probability space  define the random variables
$$\Ga^n(i,m,t)=\mmGa\bigl((z^n(i,0),0), m,t\bigr)\,.
\tag 4.5
$$
Then, following (4.3), the processes $\{z^n(t)\}$
are defined by 
$$z^n(k,t)=\inf_{i:i\le k}\lbrakk z^n(i,0)+
\Ga^n\bigl(i, k-i, t\bigr)\rbrakk\,.
\tag 4.6
$$
Our arguments use distributional bounds on the 
initial locations $z^n(i,0)$ and the variables 
$\Ga^n\bigl(i, m, t\bigr)$. In distribution 
$\Ga^n\bigl(i, m, t\bigr)$ is equal to $\mmGa(m,t)$,
so we can ignore the indices $n$ and $i$ and
 switch to $\mmGa(m,t)$ as soon as only distributional
properties are studied. 

We close this section with two comments about 
matters that came up in Section 2. 

\subsubhead 4.1 The case $q=0$ \endsubsubhead
Our results in Sections 2 and 3 are  for the case where the 
fixed equilibrium density $q$ is stricly positive. 
Here we show how the case $q=0$ reduces to 
the standard hydrodynamic setting, through a
space-time scaling of the graphical picture. 

Suppose $q=0$, and let the initial
configurations  $\eta^n(0)$ and $z^n(0)$ be 
as in (2.8) and (3.3). Define another  initial 
particle configuration by
$\ztil^n(0)=n^\beta z^n(0)$. 
Construct the  two processes $z^n(t)$ and $\ztil^n(t)$
by formula (4.3), with these Poisson processes: for
$z^n(t)$ take a realization $\Pi$ of the
rate one space-time Poisson points, 
and for $\ztil^n(t)$ use the space-time points $\Pitil$
obtained by mapping the points of $\Pi$ by
$(x,t)\mapsto(n^\beta x, n^{-\beta}t)$. By the scaling
properties of Poisson processes, $\Pitil$ is again 
a  rate one homogeneous Poisson point process. 
But notice that the difference between the evolutions 
$z^n(t)$ and $\ztil^n(t)$ is merely a stretching and 
shrinking of the space and time axes: 
$z^n(i,t)=n^{-\beta}\ztil^n(i,n^{-\beta}t)$. By 
the standard hydrodynamic result 
$n^{-1}\ztil^n([nx],nt)\to V(x,t)$, so we get that 
$n^{\beta-1}z^n([nx],n^{1+\beta}t)\to V(x,t)$. This is
 Case 2 of Theorem 4. And similarly,
Theorem 2 holds for the stick model. 

\subsubhead 4.2 The translation $[2n^\nu qt]$ \endsubsubhead
We advance here some explanation for the 
 spatial translation $[2n^\nu qt]$ in
 our theorems. 

First consider 
 the variational formula (4.3). Suppose that the 
process  is in equilibrium with $E[\eta(i,t)]=q$.
 (4.2)--(4.3) give 
$$z(0, n^\nu t)=\inf_{y\le 0}\lbrakk qy+\frac{y^2}{4n^\nu t}
+\text{[fluctuations]}\rbrakk\,.
\tag 4.7
$$
Neglecting fluctuations,
the infimum is attained at $y=-2n^\nu qt$. So
roughly speaking 
$$\aligned
z(0, n^\nu t)&=z(-[2n^\nu qt],0)+
\mmGa\bigl( (z(-[2n^\nu qt],0),0),[2n^\nu qt], n^\nu t\bigr)\\
&\qquad +\text{[fluctuations]}\,.
\endaligned
\tag 4.8
$$ 
As a sum of independents 
the term $z(-[2n^\nu qt],0)$ has fluctuations of order 
$n^{\nu/2}$, while the $\mmGa$-term has fluctuations 
of order 
$n^{\nu/3}$ (Lemma 5.2).  It is advantageous to move the translation 
$[2n^\nu qt]$ to the left-hand side of (4.8), so that we study the
dynamics of $z([2n^\nu qt], n^\nu t)$. 
Then the minimizer in (4.7) is $y=0$,
and we  get smaller fluctuations on the right-hand 
side of (4.8). 

Alternatively, we can look at the macroscopic equation
to find the right scaling and translation
 for nontrivial dynamics. 
Suppose first that $u(x,t)=q+n^{-\beta}\rho(x,t)$
satisfies the Burgers equation (2.4). Then $\rho(x,t)$
satisfies the equation
$$\rho_t+2q\rho_x+n^{-\beta}(\rho^2)_x=0\,.
$$
In the limit $n\to\infty$ this gives $\rho_t+2q\rho_x=0$
which is solved by
 a spatial translation $\rho(x,t)=\rho_0(x-2qt)$.
To get nontrivial dynamics we speed up time and set
$w(x,t)=\rho(x, n^\beta t)$ that satisfies
$$w_t+2n^{\beta}qw_x+(w^2)_x=0\,.
$$
To eliminate the $n^{\beta}$-term let $v(x,t)=w(x+2n^{\beta}qt, t)$. 
The function $v$ then solves Burgers equation 
$v_t+(v^2)_x=0$ again. Working backwards, 
$v(x,t)=\rho(x+2n^{\beta}qt,n^\beta t)$. 

To see how $v(x,t)$ should arise microscopically, 
start with the ``hydrodynamic heuristic'' $u(x,t)\approx 
(2n\e)^{-1}\sum_{|i|\le n\e}\eta([nx]+i, nt)$. From this, 
$$\rho(x,t)\approx 
\frac1{2\e}\cdot\frac1{n^{1-\beta}}\sum_{|i|\le n\e}
 \bigl\{\eta([nx]+i, nt)-q\bigr\},$$
 and then 
$$v(x,t)\approx  
\frac1{2\e}\cdot\frac1{n^{1-\beta}}\sum_{|i|\le n\e}
\bigl\{\eta([nx]+[2n^{1+\beta}qt]+i, n^{1+\beta}t)-q\bigr\}\,.$$
This is exactly what Theorem 2 states, with the 
right translation $[2n^{1+\beta}qt]$ again. 

\hbox{}

\subhead 5.\ Auxiliary lemmas \endsubhead
Throughout the proofs, we use symbols $C$, $C_1$, $C_2$, etc, 
  for constants independent
of the important  indices of the proof (such as $m$,
$n$, $i$, or $j$). 
The values of $C$, $C_1$, $C_2$, $\ldots$ may change
freely from one inequality to the next. 

We start with  an inequality for bounding the initial 
locations $z^n(i,0)$. 

\proclaim{Lemma 5.1} Suppose $\{X_i\}$ are independent
exponentially distributed random variables
with expectations $E[X_i]=q_i\in[0,b]$ where 
$b$ is a finite constant. Then  for all  $\e\in(0,1/2)$
there is a finite
constant $C=C(b,\e)>0$ such that 
for large enough $m\in\mmN$,
$$P\biggl\{ \biggl| \sum_{i=1}^m X_i- \sum_{i=1}^m q_i\biggr|
\ge \e m^{1/2+\e}\biggr\}\le 2\exp(-Cm^{2\e})\,.
\tag 5.1
$$
\endproclaim

\demo{Proof} The standard exponential Chebyshev argument.
Let $t\in(0,1/b)$. 
$$\aligned
&P\biggl\{ \sum_{i=1}^m X_i\ge  \sum_{i=1}^m q_i+
 \e m^{1/2+\e}\biggr\}\\
&\le \exp\biggl\{ -t\sum_{i=1}^m q_i-t \e m^{1/2+\e}\biggr\}
\prod_{i=1}^m E\bigl[ e^{tX_i}\bigr]\\
&=\exp\biggl\{ -t\sum_{i=1}^m q_i-t \e m^{1/2+\e}
-\sum_{i=1}^m \log(1-tq_i)\biggr\}\\
&=\exp\biggl\{ -t \e m^{1/2+\e}
+\sum_{i=1}^m \biggl( \frac{q_i^2t^2}{2} +O(t^3)\biggr)\biggr\}\\
&\le \exp\bigl\{ -t \e m^{1/2+\e}+mb^2t^2/2+O(t^3m)\bigr\}\\
&\qquad\qquad\qquad\qquad\qquad\qquad
\text{[ choose $t=b^{-2}\e m^{\e-1/2}$ ]}\\
&=  \exp\bigl\{ -\e^2b^{-2}m^{2\e}/2 +O(m^{3\e-1/2})\bigr\}\\
&\le \exp(-Cm^{2\e})\,.
\endaligned
$$
The last step is valid for $\e<1/2$. In the second equality
we expanded $\log(1-tq_i)=-tq_i-q_i^2t^2/2+O(t^3)$ where
the $O$-term is uniform over $i$ because of the uniform 
bound $q_i\le b$. The expansion is valid
because $t=b^{-2}\e m^{\e-1/2}$ can be made arbitrarily small
by restricting $m$ to be large.  The same argument 
with some sign  changes proves also the   other inequality. 
\qed
\enddemo

Next bounds on the fluctuations of the increasing paths.

\proclaim{Lemma 5.2} Suppose $a$, $s$ and $h$
are positive real numbers.

{\rm (a)} For $x\ge 2$, define 
$$I(x)=2x\cosh^{-1}(x/2)-2\sqrt{x^2-4\,}\,.
\tag 5.2
$$
Then for $a\le hs< a^2/4$, 
$$P\lbrakk \mmGa([a],s)\le \frac{a^2}{4s}-h\rbrakk \le 
\exp\lbrakk -\frac12\sqrt{a^2-4hs\,}\,I\biggl(
2+\frac{hs}{a^2}\biggr)\rbrakk\,.
\tag 5.3
$$
When $x=hs/a^2$ is small, we can use the expansion
$$I(2+x)\ge Cx^{3/2}.
\tag 5.4
$$

{\rm (b)}  There are fixed positive constants
$B_0$, $B_1$, $d_0$, $C_0$ and $C_1$ such that if
$a\ge B_0$ and 
$B_1 a^{4/3}\le hs\le d_0a^2$, then 
$$P\lbrakk \mmGa([a],s)> \frac{a^2}{4s}+h\rbrakk \le 
C_0\exp\lbrakk -C_1 \frac{s^3h^3}{a^{4}}\rbrakk\,.
\tag 5.5
$$
\endproclaim

\demo{Proof} Part (a). The random variables $\mmL(s,s)$ are
superadditive is the sense that, for any $0<s<t$, 
$$\mmL((0,0),(s,s))+\mmL((s,s),(t,t))\le \mmL((0,0),(t,t))
\qquad \text{a.s.}
$$
It follows that there exists a function $I(x)$ such that 
$$\sup_{s>0}\frac1{s}\log P\{ \mmL(s,s)\ge sx\}=
\lim_{s\to\infty}\frac1{s}\log P\{ \mmL(s,s)\ge sx\}=I(x)\,.
\tag 5.6
$$
Since $s^{-1} \mmL(s,s)\to 2$
as $s\to\infty$, $I(x)=0$ for
$x<2$. 
Kim (1996)  proved that $I(x)$ is bounded below by
the expression in (5.2), and Sepp\"al\"ainen (1998b) that
this expression equals $I(x)$. By
 (5.6) and  the observation that $\mmL(a,b)\overset{d}\to=
\mmL\bigl(\sqrt{ab\,},\sqrt{ab\,}\,\bigr)$, 
$$\aligned
&P\lbrakk \mmGa([a],s)\le \frac{a^2}{4s}-h\rbrakk
=P\lbrakk \mmL\biggl(\frac{a^2}{4s}-h,s\biggr)\ge [a]\rbrakk\\
&\le \exp\lbrakk -\frac12\sqrt{a^2-4hs\,}\,I\biggl(
\frac{2[a]}{\sqrt{a^2-4hs\,}}\biggr)\rbrakk\,.
\endaligned
$$
The argument of $I(\cdot)$ is estimated below by
$$\aligned
\frac{2[a]}{\sqrt{a^2-4hs\,}}
\ge 
\frac{2(1-1/a)}{\sqrt{1-4hsa^{-2}\,}}\ge 2\biggl(1-\frac1a\biggr)
\biggl(1+\frac{hs}{a^{2}}\biggr)\ge 2+\frac{hs}{a^{2}}\,,
\endaligned
$$
provided $hs\ge a$. 

Part (b) is a consequence of case 4
of Lemma 7.1 in Baik et al.\ (1999). We 
check the assumptions of that Lemma. 
First we express the probability (5.5) in terms of $\mmL$, then convert it 
to the $\phi_n(\lambda)$-notation of 
Baik et al.\ (1999):
$$\aligned
&P\lbrakk \mmGa([a],s)> \frac{a^2}{4s}+h\rbrakk
=P\lbrakk \mmL\biggl(\frac{a^2}{4s}+h,s\biggr)< [a]\rbrakk\\
&=P\lbrakk \mmL\biggl(\frac{a^2}{4s}+h,s\biggr)\le [a]-1\rbrakk
=\phi_{[a]-1}\biggl(\frac{a^2}{4}+hs\biggr)\,.
\endaligned
$$
According to case 4 of Lemma 7.1 in Baik et al.,
we can bound 
$$\phi_{[a]-1}\biggl(\frac{a^2}{4}+hs\biggr)\le C_0\exp(C_1t^3)
\tag 5.7
$$
with $t$ defined by the equation
$$1-\frac{t}{2^{1/3}[a]^{2/3}}=\frac{\sqrt{a^2+4hs\,}}{[a]}\,,
\tag 5.8
$$
provided 
$$
1+\frac{M_7}{2^{1/3}[a]^{2/3}}
\le
\frac{\sqrt{a^2+4hs\,}}{[a]} 
\le 1+\delta_6\,,
\tag 5.9
$$
where $M_7$ and $\delta_6$ are certain positive constants
that appear in the development of Baik et al. 

The first inequality in (5.9) is equivalent to 
$$
M_7\le 2^{1/3}\biggl( \frac{\sqrt{a^2+4hs\,}}{[a]^{1/3}}-[a]^{2/3}
\biggr)\,.
\tag 5.10
$$
Provided $hs\le 2a^2$,
the right-hand side above is bounded below by 
$$\aligned
2^{1/3}\biggl( \frac{a+hs/a}{a^{1/3}}-a^{2/3}
\biggr)
= 2^{1/3}\frac{hs}{a^{4/3}}\,.
\endaligned
\tag 5.11
$$
Thus the first inequality in (5.9) is satisfied 
if $ B_1 a^{4/3}\le hs\le 2a^2$ for a large enough constant $B_1$. 

For the second inequality in (5.9) observe that 
$$
\frac{\sqrt{a^2+4hs\,}}{[a]}\le \frac{a}{[a]}\biggl(1+\frac{4hs}{a^2}\biggr)
$$
which is $\le 1+\delta_6$ provided $a$ is large enough 
and $hs\le d_0a^2$ for a small enough $d_0$. 

We have verified the conditions of case 4
of Lemma 7.1 in Baik et al. Again because
(5.11) is below the right-hand side of (5.10), 
we see that $t$ defined by (5.8) satisfies
$-t\ge 2^{1/3}{hs}a^{-4/3}$, so inequality (5.7)
becomes (5.5).
\qed
\enddemo

\hbox{}

\demo{Remark 5.1} In our typical application of Lemma 5.2,  
$a$ and $s$ are of the same large order $m$, and 
$h$ is of the order $m^{1/3+\e}$. Then the bound
in (5.3) is $C_1\exp(-C_2m^{3\e/2})$ and in (5.5) 
 $C_1\exp(-C_2m^{3\e})$. 
\enddemo

\demo{Remark 5.2} For the growth model associated with 
the exclusion and zero-range processes, the lower tail
estimate 
(5.3) is available in Sepp\"al\"ainen (1998c) and in
Johansson (1999). But the upper tail
estimate (5.5) has not been derived at the time of writing this
paper. 
\enddemo

\hbox{}

\proclaim{Lemma 5.3} For $\e>0$, $a>0$, $\beta>0$, and $\nu\ge 1$, 
define a  deterministic quantity $R_n$
  by 
$$R_n=n^\nu tq^2+nxq+n^{1-\beta}\vinf |x|+a(n^{\nu/3+\e}+n^{1/2+\e})\,.
\tag 5.12
$$ 
Then there are finite constants $C_i>0$ such that  
$$P\bigl\{ z^n([nx]+[2n^\nu qt],n^\nu t)> R_n\bigr\} \le C_1\exp(-C_2n^{2\e})
\tag 5.13
$$ 
for all $n$. 
\endproclaim

\demo{Proof} By the variational formula (4.6) and  
Lemmas 5.1 and 5.2, 
$$\aligned
&P\bigl\{ z^n([nx]+[2n^\nu qt],n^\nu t)> R_n\bigr\}\\
&\le P\bigl\{ z^n([nx],0)\ge  nxq+n^{1-\beta}\vinf |x|
 +an^{1/2+\e} \bigr\}\\
&\qquad +P\bigl\{ \Ga^n([nx],[2n^\nu qt],n^\nu t) \ge 
n^\nu tq^2 +an^{\nu/3+\e}\bigr\}\\
&\le C_1\exp(-C_2n^{2\e})+C_1\exp(-C_2n^{3\e})\,.\qquad
\qed
\endaligned
$$
\enddemo

The main lemma of this section reduces the range of
indices that need to be considered in the variational
formula (4.6). 

\proclaim{Lemma 5.4} Suppose $\beta>0$ and $\nu\ge 1$,
and let $\xi$ be any number that satisfies $\xi\ge \nu-\beta$
and $\xi>2\nu/3$. 
Then if $b>0$ is large enough, the following holds
with probability one: for large enough $n$, 
$$\multline
z^n([nx]+[2n^\nu qt],n^\nu t)\\
=\min\lbrak
z^n(i,0)+ \Ga^n(i,[nx]+[2n^\nu qt]-i,n^\nu t)\,:\, |i-[nx]|\le [bn^{\xi}]
\rbrak.
\endmultline
\tag 5.14
$$
Furthermore, 
$$\sum_n P\lbrak\text{ {\rm (5.14) fails for $n$} } \rbrak<\infty\,.
\tag 5.15
$$
\endproclaim

\demo{Remark 5.3} The distinction between the weak inequality
$\xi\ge \nu-\beta$
and the strict inequality $\xi>2\nu/3$ in the hypothesis
 is significant. It arises
in the proof of this Lemma 5.4 and influences the error terms 
of our theorems. 
\enddemo

\demo{Proof of Lemma 5.4} We shall prove (5.14). The reader
can accumulate the estimates as we proceed and observe that
(5.15) is also true.

By Lemma 5.3 and Borel-Cantelli,   we may suppose that 
$$z^n([nx]+[2n^\nu qt],n^\nu t)\le R_n,
$$ 
at the expense of discarding
an event of probability zero and by taking $n$ large enough.
Fix $\delta_0\in(0,1)$ and take $\e$ small enough
in the definition (5.12) of $R_n$ so that 
$\nu/3+\e<\nu$ and $1/2+\e<1$. Then, since
 $\bigl|z^n([(1+\delta_0)n^\nu tq]+[nx] ,0)\bigr| $ is a sum of
$\bigl|[(1+\delta_0)n^\nu tq]+[nx]\bigr|$ independent exponential
random variables with expectations in the range
 $q\pm n^{-\beta}\vinf $, 
basic large deviation estimates show
that 
$$\sum_n P\bigl\{ z^n([(1+\delta_0)n^\nu tq]+[nx] ,0)\le R_n\bigr\}
 < \infty\,.
\tag 5.16
$$ 
Consequently we may also
assume 
$$ z^n([(1+\delta_0)n^\nu tq]+[nx] ,0)> R_n\,.
\tag 5.17
$$ 
It then follows from the variational formula (4.6) that,
almost surely, 
$$\multline
z^n([nx]+[2n^\nu qt],n^\nu t)\\
=\inf\lbrak
z^n(j,0)+ \Ga^n(j,[nx]+[2n^\nu qt]-j,n^\nu t)\,:\\
 j< [(1+\delta_0)n^\nu tq]+[nx]
\rbrak
\endmultline
\tag 5.18
$$
for large enough $n$. 

To conclude the proof we shall show that indices $j$
outside the range $|j-[nx]|\le [bn^{\xi}]$ cannot give the
infimum in (5.18). This will follow from showing that,
almost surely,  
$$\aligned
&z^n([nx]+[bn^\xi]+i ,0)+ 
\Ga^n([nx]+[bn^\xi]+i,[2n^\nu qt]-[bn^\xi]-i,n^\nu t)\\
&\ge z^n([nx],0)+ \Ga^n([nx],[2n^\nu qt],n^\nu t)
\endaligned 
\tag 5.19
$$
for all $0\le i\le [(1+\delta_0)n^\nu tq]$, and that
$$\aligned
&z^n([nx]-[bn^\xi]-i ,0)+ 
\Ga^n([nx]-[bn^\xi]-i,[2n^\nu qt]+[bn^\xi]+i,n^\nu t)\\
&\ge z^n([nx],0)+ \Ga^n([nx],[2n^\nu qt],n^\nu t)
\endaligned
\tag 5.20
$$
for all $i\ge 0$. The upper bound 
$j< 
[(1+\delta_0)n^\nu tq]+[nx]$ in (5.18) permitted us to
restrict the range of $i$ to $i\le [(1+\delta_0)n^\nu tq]$ in (5.19). 
The benefit is that the argument
 $[2n^\nu qt]-[bn^\xi]-i$ of $\Ga^n$ in (5.19) is of order $n^\nu$
throughout the range of $i$, which makes the estimation
 easier because there is no need for separate 
arguments for values of smaller order.

We estimate various terms separately in three sublemmas. 

\proclaim{Sublemma 5.1} For any fixed $b>0$ and $\delta>0$, the following
statements hold almost surely: for all large enough $n$, 
$$ \multline
\Ga^n([nx]+[bn^\xi]+i,[2n^\nu qt]-[bn^\xi]-i,n^\nu t)\\
>
\frac1{4n^\nu t}\bigl(2n^\nu qt-bn^\xi-i\bigr)^2-
\delta n^{\nu/3+\delta\nu}
\endmultline 
\tag 5.21
$$
for all $0\le i\le [(1+\delta_0)n^\nu tq]$, and 
$$  \multline
\Ga^n([nx]-[bn^\xi]-i,[2n^\nu qt]+[bn^\xi]+i,n^\nu t)\\
>
\frac1{4n^\nu t}\bigl(2n^\nu qt+bn^\xi+i\bigr)^2-
\delta\bigl( n^{\nu/3+\delta\nu}+i^{1/3+\delta} \bigr)
\endmultline 
\tag 5.22
$$
for all $i\ge 0$.
\endproclaim

\demo{Proof of Sublemma 5.1} We shall prove (5.22) and leave (5.21)
to the reader.  Their  difference is 
that in (5.22), due to the unbounded range of $i$, we need 
 $i$ explicitly in the estimates and sum over $i\ge 0$
in the end. (5.21) is easier because one estimate uniformly
over $i$ is sufficient. 

Let $A_n$ denote the event that 
(5.22) holds for all $i\ge 0$. 
Our goal is to prove
$$\sum_n P(A_n^c)<\infty
\tag 5.23
$$
so that by Borel-Cantelli $A_n$ happens for all large 
enough $n$, a.s. 
For fixed $n$ and $i$, 
the event that (5.22) fails has the same probability
as the event
$$ \mmGa([2n^\nu qt]+[bn^\xi]+i,n^\nu t)\le 
\frac1{4n^\nu t}\bigl(2n^\nu qt+bn^\xi+i\bigr)^2-
\delta\bigl( n^{\nu/3+\delta\nu}+i^{1/3+\delta} \bigr)\,.
\tag 5.24
$$
By shrinking $b$ slightly we can discard the integer
parts. We bound  the probability of
the event (5.24)  by (5.3). Now
$$\aligned
\sqrt{a^2-4hs\,}&=\biggl[
\bigl(2n^\nu qt+bn^\xi+i\bigr)^2 - 4\delta t
\bigl( n^{4\nu/3+\delta\nu}+i^{1/3+\delta}n^\nu \bigr)\biggr]^{1/2}\\
&\ge C_1(n^\nu+i)
\endaligned 
$$
and 
$$\frac{hs}{a^2}\ge C_2\,\frac{n^{4\nu/3+\delta\nu}+i^{1/3+\delta}n^\nu}
{n^\nu+i} \ge C_3n^{-2\nu/3+\delta\nu}\,, 
$$
uniformly over $i\ge 0$,  for all large enough $n$.  
Thus the probability of the event (5.24) is at most
$$\exp\bigl[ -C(n^\nu+i)I(2+C_1n^{-2\nu/3+\delta\nu})\bigr]
\le \exp\bigl( -Cn^{3\delta\nu/2}-Cin^{-\nu+3\delta\nu/2}\bigr)\,,
$$
where we applied the expansion (5.4).
Summing this over $i\ge 0$, we get that 
$$P(A_n^c)\le \sum_{i\ge 0}
\exp\bigl( -Cn^{3\delta\nu/2}-Cin^{-\nu+3\delta\nu/2}\bigr) 
\le C_1 n^{\nu(1-3\delta/2)}\exp\bigl( -Cn^{3\delta\nu/2}\bigr)\,.
$$
This last expression is summable over $n$, so (5.23) happens. 
\qed
\enddemo

\proclaim{Sublemma 5.2} For any fixed $\delta>0$, this 
 holds almost surely: for all large enough $n$, 
$$ \Ga^n([nx],[2n^\nu qt],n^\nu t)\le 
n^\nu tq^2+
\delta n^{\nu/3+\delta\nu}\,.
\tag 5.25
$$
\endproclaim

\demo{Proof of Sublemma 5.2} Deviation bound (5.5) and Borel-Cantelli. 
\qed
\enddemo

\proclaim{Sublemma 5.3} For any fixed $b>0$ and $\delta>0$, the following
statements hold almost surely: for all large enough $n$, 
$$\multline
 z^n([nx],0)-z^n([nx]+[bn^\xi]+i,0)\\
\le -\bigl( bn^\xi+i\bigr)
\bigl( q-n^{-\beta}\vinf\bigr)+
\delta\bigl( n^{\xi/2+\delta}+i^{1/2+\delta} \bigr)
\endmultline
\tag 5.26
$$
for all $0\le i\le [(1+\delta_0)n^\nu tq]$, and 
$$\multline
 z^n([nx],0)-z^n([nx]-[bn^\xi]-i,0)\\
\le \bigl( bn^\xi+i\bigr)
\bigl( q+n^{-\beta}\vinf\bigr)+
\delta\bigl( n^{\xi/2+\delta}+i^{1/2+\delta} \bigr)
\endmultline
\tag 5.27
$$
for all $i\ge 0$.
\endproclaim

\demo{Proof of Sublemma 5.3} This lemma 
is a consequence of Borel-Cantelli and the distribution
of the increments $z^n(j+1,0)-z^n(j,0)$. We prove
(5.27). The argument for (5.26) is the same. 

We can  
realize the initial locations $z^n(j,0)$ so that 
the inequalities 
$$(q-n^{-\beta}\vinf)X_j \le z^n(j,0)-z^n(j-1,0) \le
 (q+n^{-\beta}\vinf)X_j
\tag 5.28
$$
are valid, where $\{X_j\}$ are
i.i.d.\ exponential random variables with expectation 
$E[X_j]=1$. For fixed $n$ and $i$, the probability that 
(5.27) fails is
bounded above by 
$$\aligned
&P\lbrakk \sum_{j=1}^{[bn^\xi]+i} (q+n^{-\beta}\vinf)X_j \\
&\qquad\qquad\qquad\qquad > 
 \bigl( bn^\xi+i\bigr)
\bigl( q+n^{-\beta}\vinf\bigr)+
\delta\bigl( n^{\xi/2+\delta}+i^{1/2+\delta} \bigr)\rbrakk\\
&\le P\lbrakk \sum_{j=1}^{[bn^\xi]+i} X_j  > 
 \bigl( bn^\xi+i\bigr)+
\delta\frac{ n^{\xi/2+\delta}+i^{1/2+\delta}}
{ q+n^{-\beta}\vinf}\rbrakk\,.
\endaligned
\tag 5.29$$
This probability is bounded above by 
$$\exp\lbrakk - \bigl( [bn^\xi]+i\bigr)
\kappa\biggl( 1+\frac{ \delta(n^{\xi/2+\delta}+i^{1/2+\delta})}
{ (q+n^{-\beta}\vinf) ( bn^\xi+i) } \biggr)\rbrakk
\tag 5.30
$$
where $\kappa(x)=x-1-\log x$ is the Cram\'er rate function 
for the Exp(1)-distribution. In case the reader is used 
to thinking of large deviation rate functions only asymptotically
valid, note that inequality
$P\bigl(\sum_1^m X_j\ge ma\bigr)\le \exp\{-m\kappa(a)\}$ for $a>1$, and its 
lower tail counterpart, are valid for finite $m$ due to the 
supermultiplicativity
$P\bigl(\sum_1^{l+m} X_j\ge (l+m)a\bigr))\ge 
P\bigl(\sum_1^l X_j\ge la\bigr))\cdot
P\bigl(\sum_1^m X_j\ge ma\bigr))$. 

For small $x$ we have the quadratic lower bound $\kappa(1+x)\ge Cx^2$,
so (5.30) is further bounded above by
$$\exp\lbrakk -C\frac{n^{\xi+2\delta}+i^{1+2\delta}}
{bn^\xi+i} \rbrakk\,.
\tag 5.31
$$
Summing the quantities in (5.31) over $i\ge 0$, we bound the 
probability that, for a fixed $n$, (5.27) fails for {\it some}
 $i\ge 0$, by
$$\aligned
&\sum_{i\ge 0}\exp\lbrakk -C\frac{n^{\xi+2\delta}+i^{1+2\delta}}
{bn^\xi+i} \rbrakk\\
&\le \sum_{0\le i\le n^\xi}
\exp( -C n^{2\delta}) +  \sum_{ i> n^\xi}
\exp( -C i^{2\delta})\,.
\endaligned
$$
The last line above is summable over $n$. Hence by Borel-Cantelli,
it is almost surely true that for large enough $n$, 
(5.27) holds for all $i\ge 0$. 
\qed
\enddemo

We return to complete the proof of Lemma 5.4. By (5.21), (5.25), and (5.26), 
inequality (5.19) will hold for large enough $n$ if we can 
show that 
$$\multline  
\frac1{4n^\nu t}\bigl(2n^\nu qt-bn^\xi-i\bigr)^2-n^\nu tq^2 -
2\delta n^{\nu/3+\delta\nu} \\
\ge 
-\bigl( bn^\xi+i\bigr)
\bigl( q-n^{-\beta}\vinf\bigr)+
\delta \bigl(n^{\xi/2+\delta}+i^{1/2+\delta}\bigr)
\endmultline
\tag 5.32
$$
holds for all $0\le i\le [(1+\delta_0)n^\nu tq]$. 
(5.32) simplifies to 
$$\multline 
\frac{b^2}{4t}n^{2\xi-\nu}
+ \frac{b}{2t}in^{\xi-\nu}
+\frac{i^2}{4n^\nu t} -2\delta n^{\nu/3+\delta\nu}\\
\ge b\vinf n^{\xi-\beta}+ \vinf in^{-\beta}
+\delta \bigl(n^{\xi/2+\delta}+i^{1/2+\delta}\bigr) \,.
\endmultline
\tag 5.33
$$
Now suppose $\xi\ge\nu-\beta$ so that $n^{2\xi-\nu}\ge n^{\xi-\beta}$
and $n^{\xi-\nu}\ge n^{-\beta}$. 
Then (5.33)  follows 
from 
$$\multline 
\biggl( \frac{b}{4t}-\vinf\biggr)bn^{2\xi-\nu}
+\biggl( \frac{b}{2t}-\vinf\biggr)in^{\xi-\nu}
+\frac{i^2}{4n^\nu t}\\
\ge \delta\bigl( 2n^{\nu/3+\delta\nu}
+n^{\xi/2+\delta}+i^{1/2+\delta}\bigr)\,.
\endmultline
$$
Now choose $b>4t\vinf$ so that the two coefficients in 
parentheses on the left-hand side are positive and large. 
The condition $\xi>2\nu/3$ is exactly what is needed to have 
$n^{2\xi-\nu}>n^{\nu/3+\delta\nu}+n^{\xi/2+\delta}$
for all large enough $n$, if $\delta>0$ is small enough. 
The $i$-term on the right-hand side is controlled by the 
observation 
$$\delta  i^{1/2+\delta}\le C\bigl( n^{2\xi-\nu}+in^{\xi-\nu}\bigr)
$$
for all $i\ge 0$, provided $n$ is large enough. 

The argument for (5.20) goes exactly the same way. 
This completes the proof of Lemma 5.4.
\qed
\enddemo

\hbox{}

\subhead 6.\ Proof of Theorem 3\endsubhead
Let $\xi$ satisfy $\xi>2\nu/3$ and $\nu-\beta\le\xi\le\nu-\e$
for some small $\e>0$. 
Let $M<\infty$ be a large finite constant, to be chosen later. 
Define the events
$$A_n=\lbrak \,\bigl| 
z^n([nx]+[2n^\nu qt], n^\nu t)-z^n([nx],0)-n^\nu tq^2\bigr|> 2Mn^{2\xi-\nu}
\rbrak\,.
$$
By Borel-Cantelli,
 Theorem 3 will follow from proving the summability
$$\sum_nP(A_n)<\infty\,.
\tag 6.1
$$
To see this, compare $2\xi-\nu$ with the exponent in the error 
of (3.5): In the case $\nu\le 3\beta$, take $\xi=2\nu/3+\delta/3$, 
so that $2\xi-\nu<\nu/3+\delta$. 
In the case $\nu> 3\beta$, we can take $\xi=\nu-\beta$ so that
$2\xi-\nu=\nu-2\beta<\nu-2\beta+\delta$. 

 Now to prove (6.1). 
By  Lemma 5.4, for large enough $n$ it is the case that
$$\aligned
&z^n([nx]+[2n^\nu qt], n^\nu t)-z^n([nx],0)-n^\nu tq^2\\
&=\inf_{|i|\le bn^\xi}\lbrak 
z^n([nx]-i,0)-z^n([nx],0)+qi\\
&\qquad\qquad+ \Ga^n([nx]-i, [2n^\nu qt]+i, n^\nu t)-n^\nu tq^2-qi\rbrak\,.
\endaligned
\tag 6.2
$$
The probability $P(A_n)$ 
is bounded above by 
$$\aligned
&\sum_{|i|\le bn^{\xi}} P\lbrak 
\,\bigl|z^n([nx]-i,0)-z^n([nx],0)+qi\bigr|\ge   Mn^{2\xi-\nu} \rbrak\\
+\;&\sum_{|i|\le bn^{\xi}} P\lbrak
\,\bigl| \mmGa([2n^\nu qt]+i, n^\nu t)-n^\nu tq^2-qi\bigr|\ge  Mn^{2\xi-\nu}
\rbrak\\
+\;&P\lbrak \text{ (6.2) fails for $n$ }\rbrak\,.
\endaligned
\tag 6.3
$$
To bound the first probability in (6.3), apply (5.28)
to get, 
 for each $i$, 
$$\aligned
&P\lbrak 
\,\bigl|z^n([nx]-i,0)-z^n([nx],0)+qi\bigr|\ge   Mn^{2\xi-\nu} \rbrak\\
&\le P\lbrakk \sum_{j=1}^{|i|}  (q-n^{-\beta}\vinf)X_j
\le q|i| -Mn^{2\xi-\nu} \rbrakk\\
&\qquad +P\lbrakk \sum_{j=1}^{|i|}  (q+n^{-\beta}\vinf)X_j
\ge q|i| +Mn^{2\xi-\nu} \rbrakk\\
&\le P\lbrakk \sum_{j=1}^{|i|}  X_j
\le |i|+C_1n^{-\beta}|i| -M_1n^{2\xi-\nu} \rbrakk\\
&\qquad +P\lbrakk \sum_{j=1}^{|i|} X_j
\ge |i| -C_1n^{-\beta}|i|+M_1n^{2\xi-\nu} \rbrakk\,.
\endaligned
$$
In the second inequality, $C_1$ is a new constant, and 
$M_1=M/(q+1)$ accounts for the effect of dividing 
$M$ by $(q\pm n^{-\beta}\vinf)$ when $n$ is large enough.

Since $n^{-\beta}|i|\le bn^{\xi-\beta}$ and $\xi\ge\nu-\beta\Longrightarrow
2\xi-\nu\ge \xi-\beta$, by choosing $M$ large enough at the 
outset we have a further constant $M_2>0$ such that
 $M_1n^{2\xi-\nu}$ $-C_1n^{-\beta}|i|$ $\ge$ $ M_2n^{2\xi-\nu} $. 
Next, apply the large deviation rate function
$\kappa$  for Exp(1)-variables 
as in (5.29)--(5.30). The new upper bound becomes
$$\aligned
&P\lbrakk \sum_{j=1}^{|i|}  X_j
\le |i|-M_2n^{2\xi-\nu} \rbrakk
+P\lbrakk \sum_{j=1}^{|i|} X_j
\ge |i| +M_2n^{2\xi-\nu} \rbrakk\\
&\le \exp\lbrak -|i|\kappa\bigl(1-M_2n^{2\xi-\nu}|i|^{-1}\bigr)\rbrak\\
&\qquad\qquad\qquad\qquad
 +\exp\lbrak -|i|\kappa\bigl(1+M_2n^{2\xi-\nu}|i|^{-1}\bigr)\rbrak\\
&\le \exp\lbrakk M_2n^{2\xi-\nu}+|i|\log
\bigl(1-M_2n^{2\xi-\nu}|i|^{-1}\bigr)\rbrakk\\
&\qquad\qquad\qquad +
\exp\lbrakk -M_2n^{2\xi-\nu}+|i|\log
\bigl(1+M_2n^{2\xi-\nu}|i|^{-1}\bigr)\rbrakk\,.
\endaligned
\tag 6.4
$$
Check that the functions $(1/x)\log(1\pm x)$ are 
maximized by taking $x>0$ as small as possible. Thus
we get an upper bound for (6.4) by replacing 
$|i|$ by its upper bound $bn^\xi$. Expanding the 
$\log$ then gives the upper bound 
$2\exp(-Cn^{3\xi-2\nu})$ for (6.4) [here the assumption
$\xi\le\nu$ becomes useful]. Tracing backwards, we
conclude that 
$$\text{[ the first sum in (6.3) ]\, } \le C_1n^\xi\exp(-Cn^{3\xi-2\nu})\,.
\tag 6.5
$$

Now to bound the second probability in (6.3). 
Again because $|i|=O(n^\xi)$, 
$$\aligned
&P\lbrak
\,\bigl| \mmGa([2n^\nu qt]+i, n^\nu t)-n^\nu tq^2-qi\bigr|\ge  Mn^{2\xi-\nu}
\rbrak\\
&\le
P\lbrakk
\,\biggl| \mmGa([2n^\nu qt]+i, n^\nu t)-\frac1{4n^\nu t}
\bigl( 2n^\nu qt+i\bigr)^2\biggr|>  M_1 n^{2\xi-\nu}
\rbrakk\,,
\endaligned
\tag 6.6
$$
for a constant $M_1>0$, 
provided $M$ was chosen large enough. Since
$|i|\le Cn^\xi\le Cn^{\nu-\e}$, Lemma 5.2 implies that the probability
in (6.6) is at most $C_1\exp(-C_2n^{3\xi-2\nu})$. 
It follows that the bound in (6.5) works also for the second
sum in (6.3). Thus the sum in (6.1) has the following bound:
$$\sum_nP(A_n)\le \sum_n C_1n^\xi\exp(-Cn^{3\xi-2\nu})
+ \sum_n P\lbrak \text{(6.2) fails for $n$}\rbrak<\infty\,.
$$
The summability is a consequence of the assumption
$\xi>2\nu/3$ and (5.15). 
(6.1) holds, and we have proved
 Theorem 3. 

\hbox{}

\subhead 7.\ Proof of Theorem 4 \endsubhead
The proofs for the different cases 
are Borel-Cantelli estimates
for the distribution of the random variable
$z^n([nx]+[2n^\nu qt], n^\nu t)$. For the 
sake of readability, we do not 
formulate explicit probability estimates and instead 
write statements of the type 
(2.10)--(2.11). Behind each a.s.\ error estimate
is a summable deviation probability, as the 
reader can verify from the arguments. 

  Case 3 can be proved  quickly from Theorem 3: 
in (3.5) replace the term $z^n([nx],0)$ by its
expectation (3.4) plus fluctuation $o(n^{1/2+\delta})$. 
We concentrate on proving Cases 1 and 2. 

\subsubhead Proof of Theorem 4, Case 1 \endsubsubhead
Assuming  $\nu>1+\beta$, the goal is to show that
$$\lim_{n\to\infty} n^{2\beta-\nu}\lbrakk
z^n([nx]+[2n^\nu qt], n^\nu t)-n^\nu tq^2-nxq\rbrakk
=V_\infty(0,t)\,.
\tag 7.1
$$
Let $y$ be a number that achieves the infimum in (3.10)
for $x=0$. Set $i=[n^{\nu-\beta}y]$ in the expression
inside the braces on the right-hand side of (5.14). 
For large $n$, we get
the upper bound
$$\aligned
&z^n([nx]+[2n^\nu qt], n^\nu t)\\
&\le z^n([n^{\nu-\beta}y],0)+\Ga^n([n^{\nu-\beta}y],
[nx]+[2n^\nu qt]-[n^{\nu-\beta}y], n^\nu t)\\
&\le n^{\nu-\beta}yq+n^{1-\beta}V_0(n^{-1}[n^{\nu-\beta}y])\\
&\qquad +n^\nu tq^2+\frac{y^2}{4t}n^{\nu-2\beta}+nxq
-n^{\nu-\beta}yq +o(n^{\nu-2\beta}) \\
&\le n^\nu tq^2 +nxq +n^{\nu-2\beta}\biggl\{
      n^{1-\nu+\beta}V_0(n^{\nu-\beta-1}y)+\frac{y^2}{4t}\biggr\}
+o(n^{\nu-2\beta})\,.
\endaligned
\tag 7.2
$$
The following steps were taken above: For small $\e>0$, 
it is almost surely true that, for large enough $n$,  
$$\aligned
&z^n([n^{\nu-\beta}y],0)\\
&\le E\lbrak z^n([n^{\nu-\beta}y],0) \rbrak 
+\e n^{(\nu-\beta)/2+\e}\\
&= [n^{\nu-\beta}y]q+n^{1-\beta}V_0(n^{-1}[n^{\nu-\beta}y])
+\e n^{(\nu-\beta)/2+\e}\\
&\le n^{\nu-\beta}yq+n^{1-\beta}V_0(n^{\nu-\beta-1}y)
+o(n^{\nu-2\beta})\,.
\endaligned
\tag 7.3
$$
The first step above is by Lemma 5.1. 
 $\nu>3\beta$ guarantees
that $n^{(\nu-\beta)/2+\e}=o(n^{\nu-2\beta})$
 if $\e>0$ is small enough.

Similarly by 
(5.5) for large enough $n$,  
$$\aligned
&\Ga^n([n^{\nu-\beta}y],
[nx]+[2n^\nu qt]-[n^{\nu-\beta}y], n^\nu t)\\
&\le \frac1{4n^\nu t}\bigl([nx]+[2n^\nu qt]-[n^{\nu-\beta}y]\bigr)^2
+ \e n^{\nu/3+\e}\\
&\le n^\nu tq^2+\frac{y^2}{4t}n^{\nu-2\beta}+nxq
-n^{\nu-\beta}yq +C_1+C_2n^{1-\beta}+ \e n^{\nu/3+\e}\,.
\endaligned
\tag 7.4
$$
The term $C_1+C_2n^{1-\beta}$ accounts for terms left
out after expanding the square and for removal of integer parts
$[\cdot]$. The assumptions
$\nu>1+\beta$ and $\nu>3\beta$ guarantee that 
$C_1+C_2n^{1-\beta}+ \e n^{\nu/3+\e}=o(n^{\nu-2\beta})$
if $\e>0$ is small enough.
 
(7.3) and (7.4) justify the validity of (7.2) for large enough $n$, 
almost surely. Now we can prove one half of (7.1): 
$$\aligned
&\limsup_{n\to\infty}n^{2\beta-\nu}\lbrakk
z^n([nx]+[2n^\nu qt], n^\nu t)-n^\nu tq^2-nxq\rbrakk\\
&\le \limsup_{n\to\infty}\biggl\{
      n^{1-\nu+\beta}V_0(n^{\nu-\beta-1}y)+\frac{y^2}{4t}\biggr\}
=V_\infty(y,0)+\frac{y^2}{4t}
=V_\infty(0,t)\,.
\endaligned
\tag 7.5
$$
The second last equality follows from 
$$\lim_{m\to\infty}m^{-1}V_0(my)=V_\infty(y,0)\,.
\tag 7.6
$$

It remains to bound the $\liminf$ in (7.1) from below. 
By the assumptions $\nu>3\beta$ and 
$\nu>1+\beta$, we can choose a number
 $\varrho$ that satisfies 
$$1-\beta<\varrho<\nu-2\beta\,\,,\,\,\varrho>\nu/3>\beta\,\,,\,\,\text{ and }
\,\varrho>(\nu-\beta)/2\,.
\tag 7.7
$$
 Define a sequence of deterministic
numbers by 
$$r_n= n^\nu tq^2+nxq+\min_{|i|\le bn^{\nu-\beta}}\lbrakk n^{1-\beta}V_0(i/n)
+\frac{i^2}{4n^\nu t}  \rbrakk -2n^\varrho\,.
$$

\proclaim{Lemma 7.1} Almost surely, the inequality 
$z^n([nx]+[2n^\nu qt], n^\nu t) \ge r_n$
holds for large enough $n$. 
\endproclaim

Before proving the lemma, let us use it to finish the 
proof of Case 1 of Theorem 4. 
$$\aligned
&\liminf_{n\to\infty}n^{2\beta-\nu}\lbrakk
z^n([nx]+[2n^\nu qt], n^\nu t)-n^\nu tq^2-nxq\rbrakk\\
&\ge
\liminf_{n\to\infty} 
\min_{|i|\le bn^{\nu-\beta}}\lbrakk n^{1-\nu+\beta}V_0(i/n)
+\frac{i^2}{4n^{2\nu-2\beta} t}\rbrakk\\
&\qquad\qquad\qquad\qquad\qquad
\text{[ change of variable $i=n^{\nu-\beta}y$ ]}\\
&\ge
\liminf_{n\to\infty} 
\inf_{|y|\le b}\lbrakk n^{1-\nu+\beta}V_0(n^{\nu-1-\beta}y)
+\frac{y^2}{4t}\rbrakk\\
&\ge V_\infty(0,t)\,.
\endaligned
\tag 7.8
$$
To check the last inequality, let 
$n_j$ be a subsequence along which the $\liminf_{n\to\infty}$ 
is realized.  For each $j$ pick 
$y_{n_j}$ that realizes the infimum, pass
to a further convergent subsequence $y_{n_j}\to y$,  
and now consider different cases: If $y_{n_j}$ stays bounded
away from zero, it follows from (7.6) that
$$\lim_{j\to\infty}\lbrakk n^{1-\nu+\beta}_jV_0(n^{\nu-1-\beta}_jy_{n_j})
+\frac{y^2_{n_j}}{4t}\rbrakk
=  V_\infty(y,0)+y^2/4t \ge V_\infty(0,t)\,.
\tag 7.9
$$ 
And if $y=0$, Lipschitz
continuity of $V_0$ gives 
$$\bigl| n_j^{1-\nu+\beta}V_0(n_j^{\nu-1-\beta}y_{n_j})
 \bigr|\le \vinf |y_{n_j}|\longrightarrow 0\,,
$$
so in this case too the limit  in (7.9) is $V_\infty(y,0)+y^2/4t=0$.

(7.5) and (7.8) together imply (7.1), 
and thereby  prove Case 1 of Theorem. Before
moving on to Case 2, we check Lemma 7.1:

\demo{Proof of Lemma 7.1}
Abbreviate 
temporarily 
$$Z_n=\min\lbrak
z^n(i,0)+ \Ga^n(i,[nx]+[2n^\nu qt]-i,n^\nu t)\,:\, |i|\le bn^{\nu-\beta}
\rbrak.
\tag 7.10
$$
The assumption $\nu>3\beta$ permits
us to set $\xi=\nu-\beta$ in  Lemma 5.4,
so $z^n([nx]+[2n^\nu qt], n^\nu t)=Z_n$
for large enough $n$. The difference between 
$|i|\le bn^{\nu-\beta}$ in (7.10) 
and $|i-[nx]|\le bn^{\nu-\beta}$ in (5.14)
is irrelevant now because $\nu-\beta>1$ and we can always 
increase $b$. To prove Lemma 7.1, we show that 
$$\text{$Z_n\ge r_n$ holds for large enough $n$.}
\tag 7.11
$$

The complementary probability $P\{Z_n< r_n\}$ is bounded
above by the sum
$$\aligned
&\sum_{|i|\le bn^{\nu-\beta}} P\lbrak 
z^n(i,0)< qi+ n^{1-\beta}V_0(i/n)-n^\varrho\rbrak\\
+\;&\sum_{|i|\le bn^{\nu-\beta}} P\lbrakk
\Ga^n(i,[nx]+[2n^\nu qt]-i,n^\nu t)\\
&\qquad\qquad\qquad\quad< 
n^\nu tq^2+nxq -qi+\frac{i^2}{4n^\nu t}  -n^\varrho
\rbrakk\,.
\endaligned
\tag 7.12
$$
In  the first sum above the term for $i=0$ vanishes
because by construction $z^n(0,0)$ $=$ $V_0(0)$ $=$ $0$ with probability one. 
We bound
the sum over $1\le i\le bn^{\nu-\beta}$ and leave the 
matching argument for negative $i$'s to the reader. 
Let  $\{X_j\}$ be as in (5.28). 
First split the sum. 
$$\aligned
&\sum_{1\le i\le bn^{\nu-\beta}} P\lbrak 
z^n(i,0)< qi+ n^{1-\beta}V_0(i/n)-n^\varrho\rbrak\\
&\le 
\sum_{1\le i\le \e_1n^{\varrho+\beta}} P\lbrakk 
\sum_{j=0}^{i-1}(q-n^{-\beta}\vinf) X_j< qi+ n^{1-\beta}V_0(i/n)-n^\varrho   
\rbrakk\\
&\qquad +
\sum_{\e_1n^{\varrho+\beta} <i\le bn^{\nu-\beta}} P\lbrak 
z^n(i,0)< qi+ n^{1-\beta}V_0(i/n)- n^\varrho\rbrak
\endaligned
\tag 7.13
$$

To the first sum in (7.13) we apply a 
 large deviation
argument as in (5.29)--(5.30).
Pick $\e_1$,  $\e_2>0$ small enough and take $n$  large enough
so that 
$$\bigl(q-n^{-\beta}\vinf\bigr)^{-1}
\bigl(qi+ n^{1-\beta}V_0(i/n)-n^\varrho\bigr)\le i-\e_2n^\varrho$$
for all $1\le i\le \e_1n^{\varrho+\beta}$.
Then 
$$\aligned
&P\lbrakk 
\sum_{j=0}^{i-1}(q-n^{-\beta}\vinf) X_j< qi+ n^{1-\beta}V_0(i/n)-n^\varrho   
\rbrakk\\
&\le P\lbrakk 
\sum_{j=0}^{i-1} X_j< i -\e_2 n^\varrho\rbrakk\le 
\exp\bigl\{-i\kappa(1-\e_2n^{\varrho}i^{-1})\bigr\}\\ 
&\le \exp\bigl(-Cn^{2\varrho}i^{-1}\bigr) 
\le \exp\bigl(-C_1n^{\varrho-\beta}\bigr)\,.
\endaligned
$$

To the last sum in (7.13) we apply  Lemma 5.1. 
 Pick $0<\e<\varrho/(\nu-\beta)-1/2$ so that 
 $\e\in(0,1/2)$. Then for this range of $i$'s 
$$\aligned
&P\lbrak 
z^n(i,0)< qi+ n^{1-\beta}V_0(i/n)-n^\varrho\rbrak\\
&\le 
 P\lbrak 
z^n(i,0)< qi+ n^{1-\beta}V_0(i/n)-\e i^{1/2+\e}\rbrak\\
&\le C_2\exp(-C_3i^{2\e}) \le C_2\exp(-C_3n^{2\e(\varrho+\beta)})\,. 
\endaligned
$$

Combining the estimates gives 
$$\aligned
\sum_{1\le i\le bn^{\nu-\beta}} P\lbrak 
z^n(i,0)< qi+ n^{1-\beta}V_0(i/n)-n^\varrho\rbrak
\le 
 C_1 n^{\nu-\beta} \exp(-C_2n^{\gamma})
\endaligned
$$
where $\gamma>0$ is a new exponent that depends on the earlier
constants. The same bound is valid for the entire first sum in (7.12).

By Lemma 5.2, the probability in the second sum in (7.12) is at most
$$\aligned
&P\lbrakk
\mmGa([nx]+[2n^\nu qt]-i,n^\nu t)< 
\frac1{4n^\nu t}\bigl(nx+2n^\nu qt-i\bigr)^2   -C_3n^\varrho
\rbrakk\\
&\le \exp(-C_4n^{(3/2)(\varrho-\nu/3)})\,.
\endaligned
\tag 7.14
$$
The constant $C_3\in(0,1)$ appeared in front of $n^\varrho$ 
to subsume the difference between 
$n^\nu tq^2+nxq -qi+{i^2}/{(4n^\nu t)}$ in (7.12) and 
$\bigl(nx+2n^\nu qt-i\bigr)^2/(4n^\nu t)$ in (7.14). 
 Combining the estimates, 
we get
$$\aligned
\sum_n P\{Z_n< r_n\}<\infty\,.
\endaligned
$$
Borel-Cantelli now gives (7.11) and completes the proof of Lemma 7.1.
\qed
\enddemo

\hbox{}

\subsubhead Proof of Theorem, Case 2 \endsubsubhead
Assuming  $\nu=1+\beta$, the goal is now
$$\lim_{n\to\infty} n^{\beta-1}\lbrakk
z^n([nx]+[2n^\nu qt], n^\nu t)-n^\nu tq^2-nxq\rbrakk
=V(x,t)\,.
\tag 7.15
$$
Let $y$ be a number that achieves the infimum in (2.13)
so that 
$V(x,t)=V_0(y)+(x-y)^2/4t$. Set $i=[ny]$ in the expression
inside the braces on the right-hand side of (5.14). Repeat
the calculation in (7.2), to  get
an upper bound for large $n$:
$$\aligned
&z^n([nx]+[2n^\nu qt], n^\nu t)\\
&\le z^n([ny],0)+\Ga^n([ny],
[nx]-[ny]+[2n^\nu qt], n^\nu t)\\
&\le nyq+n^{1-\beta}V_0(y)+ n^{1/2+\e}\\
&\qquad +n^\nu tq^2+\frac{(x-y)^2}{4t}n^{1-\beta}+nxq
-nyq + n^{(1+\beta)/3+\e} \\
&\le n^\nu tq^2 +nxq +n^{1-\beta}V(x,t)
+o(n^{1-\beta})\,.
\endaligned
\tag 7.16
$$
The above steps are justified 
by Lemmas 5.1 and 5.2, as was done in (7.3)--(7.4) below.
Again $\e>0$ needs to be small enough. The estimate 
$n^{1/2+\e}+n^{(1+\beta)/3+\e}=o(n^{1-\beta})$ follows
from $\beta<1/2$, which itself is a consequence of the 
assumptions $\nu=1+\beta$ and $\nu>3\beta$. 
This gives one half of the goal (7.15). 

For the  other half of the proof we can also
follow the argument of Case 1. 
Since $\beta<1/2$, 
we can choose $\varrho$ so that 
$$\frac{1+\beta}3<1/2<\varrho<1-\beta\,.$$
 By the variational formula (2.13), for all $i$, 
$$\aligned
&n^\nu tq^2 +nxq +n^{1-\beta}V(x,t)-2n^\varrho\\
&\le n^\nu tq^2 +nxq +n^{1-\beta}V_0(i/n)
+n^{1-\beta}\frac{(x-i/n)^2 }{4t}-2n^\varrho
\\
&=\biggl[ qi+n^{1-\beta}V_0(i/n)-n^\varrho\biggr]
+ \biggl[ \frac1{4n^\nu t}\bigl( nx+2n^\nu qt-i\bigr)^2-n^\varrho\biggr]\,.
\endaligned
$$
Now the argument of Lemma 7.1 can be repeated to 
conclude that almost surely, for large enough $n$, 
$$
z^n([nx]+[2n^\nu qt], n^\nu t)
\ge n^\nu tq^2 +nxq +n^{1-\beta}V(x,t)-2n^\varrho\,,
$$
which   together with (7.16) implies (7.15). Case 2 is proved, and thereby
Theorem 4.

\hbox{}

\hbox{}

\head References \endhead 

\hbox{}

\flushpar  
D. Aldous and P. Diaconis (1995). 
Hammersley's interacting particle process and
longest increasing subsequences. 
 \ptrf \   103 199--213.

\hbox{}

\flushpar
Esposito, R., Marra, R., and Yau, H.-T.  (1994).
Diffusive limit of asymmetric simple exclusion. 
Rev. Math. Phys. 6, no. 5A, 1233--1267.

\hbox{}

\flushpar  
L. C. Evans (1998). Partial Differential Equations. 
American Mathematical Society.

\hbox{}

\flushpar  
P. A. Ferrari and L. R. G. Fontes (1994). Current
fluctuations for the asymmetric simple exclusion
process. \ap \ 22 820--832.

\hbox{}

\flushpar  
Johansson, K. (1999). Shape fluctuations and random
matrices. Preprint, \break
 math.CO/9903134.

\hbox{}

\flushpar  
Baik, J., Deift, P. and Johansson, K. (1999). 
On the distribution of the length of the longest
increasing subsequence of random permutations. 
J. Amer. Math. Soc., to appear.

\hbox{}

\flushpar  
Kim, J. H. (1996). On increasing subsequences of 
random permutations. J. Combin. Theory Ser. A 76 148--155.

\hbox{}

\flushpar  
C. Kipnis and C. Landim (1999). Scaling Limits 
of Interacting Particle Systems.
 Grundlehren der mathematischen Wissenschaften, vol 320,
Springer Verlag, Berlin.

\hbox{}

\flushpar  
T. Sepp\"al\"ainen (1996). A microscopic
model for the Burgers equation and longest increasing
subsequences.   Electronic J.
Probab.  1,  Paper 5, 1--51.

\hbox{}

\flushpar  
T. Sepp\"al\"ainen (1998a). Hydrodynamic scaling, convex
duality, and asymptotic shapes of growth models. 
 \mprf \ 4 1--26. 

\hbox{}

\flushpar  
T. Sepp\"al\"ainen (1998b). 
Large deviations for increasing sequences on the plane. 
\ptrf \ 112 221-244. 

\hbox{}

\flushpar  
T. Sepp\"al\"ainen (1998c). 
Coupling the totally asymmetric simple exclusion process
with a moving interface.
Proceedings of I Escola Brasileira de Probabilidade 
(IMPA, Rio de Janeiro, 1997). 
 \mprf \ 4 593--628.



\hbox{}

\flushpar  
Vershik, A. M., and Kerov, S. V. (1977). 
Asymptotics of the Plancherel measure of the symmetric
group and the limiting form of Young tables. 
Soviet Math. Dokl. 18 527--531.

\enddocument

\flushpar
Department of Mathematics\hfill\break 
Iowa State University\hfill\break 
Ames, Iowa 50011\hfill\break 
USA \hfill\break
seppalai\@iastate.edu\hfill\break

\enddocument

\enddocument